%% file: crepresconj.tex
\documentclass[oneside,9pt]{amsart}

\usepackage[light,onlyrm,notext,nosf,sfmathbb,frenchstyle]{kpfonts}
\numberwithin{equation}{section}
\usepackage[mathcal]{eucal}
\DeclareSymbolFont{usualmathcal}{OMS}{cmsy}{m}{n}
\DeclareSymbolFontAlphabet{\mathocal}{usualmathcal}
\usepackage[]{dsfont}
\usepackage[margin=.9in]{geometry}
\usepackage{xcolor}
\newcommand{\myName}{John Calabrese}
\newcommand{\myTitle}{On the Crepant Resolution Conjecture for Donaldson-Thomas Invariants}
\usepackage[british]{babel}
\usepackage{microtype}
\usepackage{hyperref}
	\hypersetup{colorlinks=false, pdfauthor={\myName}, pdfcreator={\myName}, pdfsubject={math.ag}, pdftitle={\myTitle}, linkcolor=blue, citecolor=brown}
\usepackage{tikz}
		\usetikzlibrary{arrows,matrix,shapes,snakes}
\usepackage[backgroundcolor=yellow,textsize=tiny,]{todonotes}
\DeclareMathOperator{\Per}{Per}
\DeclareMathOperator{\PHilb}{P-Hilb}
\newcommand{\pepe}[1]{{^p\!{#1}}}

\newcommand{\orb}[1]{\curvy{#1}}
\newcommand{\pphont}[1]{{\scshape #1}}
\newcounter{steps}
\newcommand{\stepz}{\stepcounter{steps}\paragraph{\pphont{Step \arabic{steps}}}}
\newcommand{\perv}[2]{{^{\scriptscriptstyle #1}{#2}}}
\newcommand{\qerv}[2]{{^{\scriptscriptstyle #1}{#2}}}

\renewcommand{\D}{\mathbb{D}}
\newcommand{\twit}{}
\newcommand{\otwit}[1]{\hat{#1}}
\newcommand{\correct}[1]{\underline{#1}}

\input{macros}

\newcommand{\X}{\orb{X}}

\usepackage{tikz-cd}
\begin{document}
	\author{{\tiny \text{\myName}}}
	\title{{\Large \normalfont \scshape \myTitle}}
	\address{Mathematical Institute\\ University of Oxford\\ UK}
	\email{john.robert.calabrese@gmail.com}
	\maketitle
	\begin{abstract}
		We prove a comparison formula for curve-counting invariants in the setting of the McKay correspondence, related to the \emph{crepant resolution conjecture} for Donaldson-Thomas invariants.
		The conjecture is concerned with comparing the invariants of a (hard Lefschetz) Calabi-Yau orbifold of dimension three with those of a specific crepant resolution of its coarse moduli space.
		We prove the conjecture for point classes and give a conditional proof for general curve classes.
		We also prove a variant of the formula for curve classes.
		Along the way we identify the image of the standard heart of the orbifold under the Bridgeland-King-Reid equivalence.
	\end{abstract}
	\tableofcontents
	\vspace{-.7cm}
	
	\input{intro}

	\input{one}

	\input{two}

	
	\bibliographystyle{myamsalpha.bst}
	\bibliography{biblio.bib}
\end{document}

%% file: macros.tex
\newcommand{\into}{\hookrightarrow}

\newcommand{\onto}{\twoheadrightarrow}

\newcommand{\st}{\,\middle\vert\,}

\newcommand{\C}{\mathds{C}}
\renewcommand{\O}{\mathcal{O}}
\newcommand{\R}{\mathbb{R}} 
\newcommand{\Z}{\mathbb{Z}}

\newcommand{\cat}[1]{\mathsf{#1}}
\newcommand{\curly}[1]{\mathscr{#1}}
\newcommand{\curvy}[1]{\mathcal{#1}}

\DeclareMathOperator{\ch}{ch}
\DeclareMathOperator{\Coh}{Coh}

\DeclareMathOperator{\DT}{DT}

\DeclareMathOperator{\exc}{exc}
\DeclareMathOperator{\Ext}{Ext}
\DeclareMathOperator{\Hom}{Hom}
\DeclareMathOperator{\Hilb}{Hilb}
\DeclareMathOperator{\mr}{mr}

\DeclareMathOperator{\supp}{supp}
\DeclareMathOperator{\topp}{top}
\renewcommand{\top}{\text{top}}

\usepackage{thmtools}
\declaretheoremstyle[
spaceabove=6pt, spacebelow=6pt,
headfont=\normalfont\scshape,
notefont=\normalfont,
notebraces={(}{)},
bodyfont=\normalfont,
postheadspace=.5em,
headpunct={ --},
headformat=swapnumber
]{tstyle}
\declaretheorem[style=tstyle,name=Theorem,sibling=equation]{thmm}

\declaretheorem[numberlike=equation,style=tstyle,name=Lemma]{lem}
\declaretheorem[numberlike=equation,style=tstyle,name=Definition]{defn}
\declaretheorem[numberlike=equation,style=tstyle,name=Corollary]{cor}

\declaretheoremstyle[
spaceabove=6pt, spacebelow=6pt,
headfont=\small\normalfont\itshape,
notefont=\small\normalfont\itshap,
notebraces={(}{)},
bodyfont=\small\normalfont,
postheadspace=1em,
headpunct=:,
qed=$\blacksquare$
]{pstyle}
\declaretheorem[numbered=no,style=pstyle,name=Proof]{prf}

\declaretheorem[sibling=equation,name=Remark,style=remark]{rmk}

\declaretheorem[numbered=no,name=Remark,style=remark]{rmk*}

\declaretheorem[sibling=equation,name=Situation,thmbox=M
]{situ}

\declaretheoremstyle[
	spaceabove=6pt, spacebelow=6pt,
	headfont=\normalfont\itshape,
	notebraces={(}{)},
	bodyfont=\normalfont,
	postheadspace=.5em,
	headformat=swapnumber
]{rstyle}
\declaretheoremstyle[
	spaceabove=6pt, spacebelow=6pt,
	headfont=\normalfont\itshape,
	notebraces={(}{)},
	bodyfont=\normalfont,
	postheadspace=.5em,
	headformat=swapnumber
]{rstyle}
\declaretheorem[numbered=no,name=Notation,style=rstyle]{notation}


%% file: intro.tex
\section*{Introduction} 
\label{sec:introduction} 
This paper is concerned with the \emph{crepant resolution conjecture} for Donaldson-Thomas (DT) invariants as stated in \cite[Conjectures 1 and 2]{crc}.
Our goal is to give a full proof of Conjecture 2 and a conditional proof of Conjecture 1.
We also prove of a variant of Conjecture 1 for ``partial'' DT invariants.

The goal of these conjectures is to pin down the relationship between the DT invariants of a CY3 orbifold $\orb{X}$, satisfying the hard Lefschetz condition, and the DT invariants of a natural crepant resolution $Y \to X$ of its coarse moduli space $X$.
Concretely, we relate counting invariants of $\orb{X}$ and $Y$.
The proof employs a derived equivalence between $\orb{X}$ and $Y$, which is a ``global'' version of the McKay correspondence of Bridgeland-King-Reid \cite{bkr,chentseng}.
We prove that the image of the heart $\Coh(\orb{X})$ under this equivalence is Bridgeland's category of perverse coherent sheaves $\Per(Y/X)$ \cite{tomflops}.

Before writing down the formula in symbols, it is profitable to spend a few words on the setup of the conjecture.\footnote{The reader interested in more background on DT theory (and curve-counting in general) could start from \cite{13}.}
Given a smooth and projective Calabi-Yau\footnote{\label{sign}For us \emph{Calabi-Yau} means having trivial canonical bundle $\omega_M \cong \O_M$ and torsion fundamental group $H^1(M,\O_M)=0$.} threefold $M$, we can define the DT invariants of $M$ as weighted Euler characteristics\footnote{There is a minor sign issue in this definition. We expand upon it in Remarks \ref{sign issue crc} and \eqref{sign issue again crc}. A quick inspection will show that the main formulae we prove hold regardless of sign conventions.}
\begin{align*}
	\DT_M(\beta,n) := \chi_\text{top} \left( \Hilb_M(\beta,n), \nu \right) = \sum_{k \in \Z} k \chi_{\topp} \left( \nu^{-1}(k) \right)
\end{align*}
where $\chi_\text{top}$ is the topological Euler characteristic, $\beta \in N_1(M)$ is the homology class of a curve, $n$ is an integer, $\Hilb_M(\beta,n)$ is the Hilbert (or Quot) scheme parameterising quotients of $\O_M \onto E$ of class $(\ch_0 E, \ch_1 E, \ch_2 E,\ch_3 E) = (0,0,\beta,n)$ and $\nu$ is Behrend's microlocal function \cite{behrend}.
We formally package these numbers into a generating series.
\begin{align*}
	DT(M) := \sum_{(\beta,n) \in N_1(M) \oplus \Z} DT_M(\beta,n) q^{(\beta,n)}
\end{align*}
Taking Chern characters (and using \cite[Lemma 2.2]{tomcc}) we can we replace $N_1(M) \oplus \Z$ with the numerical Grothendieck group.
To be precise, we let $N(M)$ be the K-group of coherent sheaves on $M$ modulo numerical equivalence and we define $F_1N(M)$ to be the subgroup spanned by sheaves supported in dimension at most one.
It follows that $DT(M)$ can alternatively be indexed by $F_1N(M)$:
\begin{align*}
	DT(M) = \sum_{\alpha \in F_1N(M)} DT_M(\alpha) q^\alpha
\end{align*}
and we will switch between one indexing and the other depending on circumstances. 
There is also a subgroup $F_0N(M)$ spanned by sheaves supported in dimension zero and we can define
\begin{align*}
	DT_0(M) := \sum_{\alpha \in F_0N(M)} DT_M(\alpha) q^\alpha.
\end{align*}

Let now $\orb{X}$ be a projective Calabi-Yau orbifold of dimension three and let $X$ be its coarse moduli space.
By \cite{bkr,chentseng} there is a crepant resolution $Y \to X$ of $X$ given by an appropriate Hilbert scheme of points of $\orb{X}$.
\begin{center}
	\begin{tikzpicture}
		\matrix (m) [matrix of math nodes, row sep=3em, column sep=3em, text height=1.5ex, text depth=0.25ex]
		{
		Y & & \orb{X} \\
		& X & \\
		};
		\path[->,font=\scriptsize]
		(m-1-1) edge node[auto,swap]{$f$} (m-2-2)
		(m-1-3) edge node[auto]{$g$} (m-2-2)
		;
	\end{tikzpicture}
\end{center}
The global McKay correspondence tells us that, moreover, $Y$ and $\orb{X}$ are derived equivalent via Fourier-Mukai transforms
\begin{align*}
	\Phi: D(Y) \rightleftarrows D(\orb{X}):\! \Psi
\end{align*}
inducing isomorphisms between the corresponding (numerical) K-groups.
We also assume $f$ to have fibres of dimension at most one.
This is equivalent to requiring $\X$ to be \emph{hard Lefschetz} (cfr.~\cite[Lemma 24]{bg}).

Let $F_{\exc}N(Y) \subset F_1N(Y)$ be the subgroup spanned by sheaves whose support is contracted to a point by $f$.
We have a corresponding DT series
\begin{align*}
	DT_{\exc}(Y) := \sum_{\substack{(\beta,n) \in N_1(Y) \oplus \Z \\ f_* \beta = 0}} DT_Y(\beta,n)q^{(\beta,n)}
\end{align*}
and a variant which will be useful later.
\begin{align*}
	DT_\text{exc}^\vee(Y) := \sum_{\substack{(\beta,n) \in N_1(Y) \oplus \Z \\ f_* \beta = 0}} DT_Y(-\beta,n)q^{(\beta,n)}
\end{align*}
Over $\X$ we define\footnote{
	The subscript $_\text{mr}$ stands for multi-regular, see \cite{crc}.}
$F_{\mr}N(\X) \subset F_1N(\X)$ to be the image, under the McKay correspondence functor $\Phi$, of $F_1N(Y)$.
We draw a diagram expressing the compatibilities among the various classes
\begin{equation}\tag{$\Delta$}\label{delta}
\begin{tikzcd}
	F_0N(Y) \ar[hook]{r} &
	F_{\exc}N(Y)	\ar[-, double equal sign distance]{d}{\Phi}	\ar[hook]{r}	&
	F_1N(Y)			\ar[-, double equal sign distance]{d}{\Phi}	& \\
										&
	F_0N(\orb{X})	\ar[hook]{r} &
	F_{\mr}N(\orb{X}) \ar[hook]{r} &
	F_1N(\orb{X})
\end{tikzcd}
\end{equation}
and we write down the corresponding DT series.
\begin{align*}
	DT_\text{mr}(\orb{X}) := \sum_{\alpha \in F_\text{mr}N(\orb{X})} DT_\orb{X}(\alpha) q^\alpha,
	\quad
	\quad
	DT_0(\orb{X}) := \sum_{\alpha \in F_0N(\orb{X})} DT_\orb{X}(\alpha) q^\alpha
\end{align*}
In \cite{crc} two formulae are conjectured to hold.
\begin{align}
	\tag{C1}\label{c1}\frac{ DT_{\mr}(\X) } { DT_0(\X) } &= \frac{ DT(Y) } {DT_{\exc}(Y)}
\end{align}
\begin{align}
	\tag{C2}\label{c2}DT_0(\X) &= \frac{ DT_{\exc}(Y) DT_{\exc}^\vee(Y) } { DT_0(Y) }
\end{align}
We give here a proof of \eqref{c2}.
In the case of transverse $A$-singuarities, this identity was proved by Jim Bryan in the appendix to \cite{dt0}.
Shortly before the last version of the present paper appeared on the arxiv, a proof of \eqref{c1} in the toric case with $A$-singularities was given by Dustin Ross \cite{dusty}.

Given \eqref{c2}, we see that proving \eqref{c1} is equivalent to showing the following.
\begin{align}
	\tag{C0}\label{c0} \DT_{\mr}(\X) &= \frac{ DT(Y) DT^\vee_{\exc}(Y) }{DT_0(Y)}
\end{align}
We give a conditional proof of \eqref{c0} (and thus of \eqref{c1}).
We also provide an unconditional proof of a variant of \eqref{c0}
\begin{align}
	\tag{C0$^\partial$}\label{c0del} \DT_{\mr}^\partial(\X) &= \frac{ DT^\partial(Y) DT^\vee_{\exc}(Y) }{DT_0(Y)}
\end{align}
for invariants $DT^\partial$ defined by taking the weighted Euler characteristic of appropriate open subschemes of $\Hilb$ and $\PHilb$, see Corollary \ref{corollario}.

\subsubsection*{A sketch of the proof} The key result is identifying the image (via $\Psi$) of $\Coh (\orb{X})$  inside $D(Y)$.
It turns out that $\Psi(\Coh(\orb{X}))$ is none other than Bridgeland's heart of perverse coherent sheaves $\Per(Y/X)$.
The relationship between $\Per(Y/X)$ and DT invariants was studied in \cite{cala} (and previously in \cite{todaflops}).
What follows contains the main ideas of the proof although it glosses over an issue of signs -- this is explained in Remarks \ref{sign issue crc}, \ref{sign issue again crc}.
As $\O_Y \in \Per(Y/X)$ one has a \emph{perverse Hilbert scheme} $\PHilb_{Y/X}(\alpha)$ parameterising quotients of $\O_Y$ in $\Per(Y/X)$ of numerical class $\alpha$.
One can then define
\begin{align*}
	DT_{Y/X}(\alpha) := \chi_\text{top}\left( \PHilb_{Y/X}(\alpha), \nu \right)
	\quad \text{ and } \quad
	&DT(Y/X) := \sum_{\alpha \in F_1N(Y)} DT_{Y/X}(\alpha) q^\alpha \\
	& DT_{\exc}(Y/X) := \sum_{\alpha \in F_{\exc}N(Y)} DT_{Y/X}(\alpha) q^\alpha.
\end{align*}
The Fourier-Mukai transform $\Psi$ not only identifies $\Coh(\orb{X})$ with $\Per(Y/X)$ but also the corresponding Hilbert schemes, so that we have $\Hilb_\orb{X}(\alpha) = \PHilb_{Y/X}(\psi({\alpha}))$.
One then obtains for free an identification between the generating series:
\begin{align*}
	DT_{\mr}(\orb{X}) = DT(Y/X) \quad \text{ and } \quad DT_0(\X) = DT_{\exc}(Y/X).
\end{align*}
In \cite[Theorem 4.4]{cala} the following relation between $DT(Y/X)$ and ordinary DT invariants was proved
\begin{align}\label{flopfo}\tag{$\bigstar$}
	DT_{\exc}(Y/X) = \frac{DT_\text{exc}^\vee (Y) DT_{\exc}(Y)}{\DT_0(Y)}
\end{align}
which instantly implies \eqref{c2}, i.e.~\cite[Conjecture 2]{crc}.

To establish \cite[Conjecture 1]{crc} (in other words, to show \eqref{c0}) it would suffice to have the formula
\begin{align}
	\tag{$\bigstar\bigstar$}\label{c8} DT(Y/X) = \frac{DT_\text{exc}^\vee (Y) DT(Y)}{\DT_0(Y)}
\end{align}
which, however, is not proved in \cite{cala}.
\begin{rmk*}
	The reader familiar with work of Toda, might notice that \eqref{c8} is implied by \cite[Theorem 7.3]{todaflops}.
	Unfortunately, the author has confirmed (through private communication during the 2014 GRIFGA/Lebesgue summer school held in Nantes) that the formula appearing in that theorem should be modified.
	In other words, the techniques of \cite{todaflops} (just as those of \cite{cala}) do \emph{not} in fact extend verbatim to our present setting.
\end{rmk*}

However, a variant of \eqref{c8} was proved \cite[Theorem 3.30]{cala}.
To state it, we need to make a few definitions.
Let $\PHilb^\partial_{Y/X}(\alpha)$ denote the open subspace of $\PHilb_{Y/X}(\alpha)$, parameterising epimorphisms $\O_Y \to E$ with $\dim \supp E \leq 1$.
Recall that the support of a complex is defined to be the union of the supports of its cohomology sheaves.
Since $X$ is allowed to have one-dimensional singular locus, one cannot detect whether $E \in \Per(Y/X)$ is supported on a curve simply by looking at its numerical class.
By taking the weighted Euler characteristic of $\PHilb^\partial$ we obtain ``partial'' DT invariants $DT^\partial(Y/X)$.
On the orbifold side, let $\Hilb_{\X}^\partial$ be the image of $\PHilb_{Y/X}^\partial$ under $\Phi$, together with the corresponding partial invariants on $\X$, $DT^\partial(\X)$.
It tautologically follows that
$$ DT^\partial_{\mr}(\X) = DT^\partial(Y/X) $$
with the usual identification of numerical classes under the McKay correspondence.
Let now $\Hilb_Y^\partial \subset \Hilb_Y$ be the open subset parameterising those quotients $\O_Y \onto E$, such that the following condition holds: if we view $\O_Y \to E$ as a morphism in $\Per(Y/X)$, then the perverse cokernel (i.e.~the cokernel in the abelian category $\Per(Y/X)$) lies in $\Coh_{\leq 1}(Y)[1]$.
In general, the perverse cokernel of $\O_Y \to E$ is a sheaf shifted by $1$; here we add the condition that it is the shift of a sheaf supported in dimension at most one (see \cite[Lemma 1.4]{cala}).
Let us denote the invariants we obtain by $DT^\partial(Y)$.
The ``partial'' version of \eqref{c8} holds:
\begin{align*}
	DT^\partial(Y/X) &= \frac{ DT^\partial(Y) DT^\vee_{\exc}(Y) } {DT_0(Y)}
\end{align*}
from which the following variants of \eqref{c0} and \eqref{c1} \cite[Conjecture 1]{crc}
\begin{align*}
	DT^\partial_{\mr}(\X) = \frac{ DT^\partial(Y) DT^\vee_{\exc}(Y) } {DT_0(Y)},
	\quad \quad \quad
	\frac { DT^\partial_{\mr}(\X) }{DT_0(\X)} = 
	\frac { DT^\partial(Y) }{ DT_{\exc}(Y) }
\end{align*}
are deduced.
Although not relevent for the present paper, we point out that when $X$ has zero-dimensional singular locus then all the ``partial'' moduli spaces coincide with the ordinary ones and thus $DT^\partial = DT$ throughout.

\bigskip
We also mention a different approach to prove the crepant resolution conjecture proposed by Bryan-Steinberg in \cite{david}.
They develop new invariants, which are a relative version of the \emph{stable pair} invariants of Pandharipande-Thomas (PT), and prove a DT/PT comparison formula.
The goal would be to relate the Bryan-Steinberg invariants of $Y/X$ with the PT invariants of $\X$ and then use the DT/PT formula for $\X$ announced by Arend Bayer.
Unfortunately, a direct comparison using the Fourier-Mukai transform $\Phi$ does not seem to work.
We feel this issue must be related to the difference between $\DT$ and $\DT^\partial$.

\subsubsection*{Acknowledgements} The author would like to thank Tom Bridgeland for precious help, Dominic Joyce for useful conversations, Arend Bayer for sharing a preliminary version of an upcoming paper and Jim Bryan and David Steinberg for sharing their ideas.

\subsubsection*{Structure of the paper} The paper is divided into two sections.
The first one is the core, as it contains the proof of the fact that $\Coh (\orb{X})$ is sent to $\Per(Y/X)$ via the derived equivalence.
In the second section we apply this result to Donaldson-Thomas invariants.

\subsubsection*{Conventions} We work over the field of complex numbers $\C$. For a scheme (or stack) $M$, $D(M)$ will denote the \emph{bounded} derived category of coherent $\O_M$-modules.


%% file: one.tex
\section{The Equivalence between $\Per(Y/X)$ and $\Coh (\orb{X})$} 
\label{sec:the_equivalence_between_per_y_x_and_coh_x_}
We work in the following setup.
\begin{situ}\label{situ1}
	Let $\orb{X}$ be a smooth, quasi-projective, Deligne-Mumford stack of dimension $n$.
	Assume the canonical bundle $\omega_{\orb{X}}$ to be Zariski-locally trivial and denote by $X$ the coarse moduli space of $\orb{X}$.
\end{situ}
\begin{rmk}
	The bundle $\omega_\orb{X}$ on $\orb{X}$ is \emph{Zariski}-locally trivial if there exists a Zariski open cover $\orb{X}' \to \orb{X}$ (where we allow $\orb{X}'$ to be a \emph{stack}) such that the restriction $\omega_\orb{X}\vert \orb{X}'$ is trivial.
	This is a technical condition which, by working locally on the coarse space $X$, allows us to reduce to the setting of \cite{bkr}.
	In fact, in the case where $\orb{X} = [V/G]$, it amounts to requiring that the canonical bundle of $V$ be $G$-equivariantly locally trivial.
	This condition seems to be missing in \cite{chentseng}.
\end{rmk}

It is beneficial to recall the framework of \cite{chentseng}.
A candidate for a resolution of $X$ (and a replacement for the equivariant Hilbert scheme found in \cite{bkr}) is given by the irreducible component $Y$ of the Hilbert scheme $\Hilb (\orb{X})$ containing the non-stacky points of $\orb{X}$.\footnote{It is probably helpful to remark that for a stack $\orb{X}$ there might be some ambiguity in the term \emph{Hilbert scheme} (see \cite{rydh}). However, we shall always interpret Hilbert schemes as Quot functors, which for Deligne-Mumford stacks were studied by Olsson and Starr \cite{os}.}
The morphism $g\colon \orb{X} \to X$ induces a morphism $\Hilb(\orb{X}) \to \Hilb(X)$ and, by restriction, a morphism $f\colon Y \to X$. We draw a diagram.
\begin{center}
	\begin{tikzpicture}
		\matrix (m) [matrix of math nodes, row sep=3em, column sep=3em, text height=1.5ex, text depth=0.25ex]
		{
		& Y \times \orb{X} &\\
		Y & &\orb{X}\\
		& X &\\
		};
		\path[->,font=\scriptsize]
		(m-2-1) edge node[auto,swap]{$f$} (m-3-2)
		(m-2-3) edge node[auto]{$g$} (m-3-2)
		(m-1-2) edge node[auto,swap]{$\pi_Y$} (m-2-1)
				edge node[auto]{$\pi_\orb{X}$} (m-2-3)
		;
	\end{tikzpicture}
\end{center}
	Under the additional assumption that $Y\times_X Y$ is at most of dimension $n+1$ it is proved in \cite{chentseng} that $Y$ is smooth and that $f$ is a crepant resolution.
	Furthermore, the scheme $Y$ represents a moduli functor and its corresponding universal object is a quotient $\O_{Y \times \orb{X}} \onto \O_\orb{Z}$.
	Finally, it is shown that one has a Fourier-Mukai equivalence $D(Y) \simeq D(\orb{X})$ with kernel given by $\O_\orb{Z}$.
	
	We recall three key results involved in the proof: the Hilbert scheme $\Hilb (\orb{X})$ commutes with \'etale base-change on $X$ \cite[Proposition 2.3]{chentseng}; \'etale-locally on $X$ the space $\orb{X}$ is isomorphic to a quotient stack $[V/G]$, with $V$ smooth and affine and $G$ a finite group (whose coarse space is thus the quotient $V/G$) \cite[Lemma 2.2.3]{vista}; the Hilbert scheme of $[V/G]$ is isomorphic to Nakamura's $G$-equivariant Hilbert scheme $G\text{-}\Hilb (V)$ \cite[Lemma 2.2]{chentseng}.
	Using these facts one reduces to \cite{bkr}, as checking that the given kernel produces an equivalence may be done locally \cite[Proposition 3.3]{chentseng}.
\begin{rmk}\label{otherperversity}
	As is usual with integral transforms, the kernel $\O_\orb{Z}$ may be interpreted as giving a functor in two different directions.
	The standard Mukai-implies-McKay convention is to take $\O_\orb{Z}$ to define a functor $\twit{\Phi}\colon D(Y) \to D(\orb{X})$ \cite{bkr,chentseng}.
	To deal with a technical issue (caused by \cite{cala}), we will also consider $\hat{\Phi} = \D\Phi\D \colon D(Y) \to D(\orb{X})$, where $\D = \R\underline{\Hom}(-,\O)$ is the duality functor.
	We denote by $\Psi$ the inverse of $\Phi$ and by $\hat{\Psi}$ the inverse of $\hat{\Phi}$.
	When $Y$ and $X$ are projective, the relationship between $\Phi$ and $\hat{\Phi}$ is quite simple, as $\hat{\Psi}$ is given by the Fourier-Mukai transform with kernel $\O_\orb{Z}$ (this is a standard consequence of \cite[Propositions 1.13 and 1.15]{nahm}).
\end{rmk}
We now briefly remind the reader of Bridgeland's heart of perverse coherent sheaves \cite{tomflops}.
In some sense, it is a reflection of the ambiguity revolving around the kernel $\O_\orb{Z}$ that we consider both the $-1$ and $0$ perversity.
The category $\perv{p}{\Per(Y/X)}$ of \emph{perverse coherent} of \emph{perversity} $p \in \{-1,0\}$ consists of those complexes $E \in D(Y)$ satisfying
\begin{itemize}
	\item $\R f_* E \in \Coh (X)$,
	\item $\Ext^{-i}_Y (E, C) = 0 = \Ext_Y^{-i}(C,E)$, for all $i > p$ and all $C \in \Coh (Y)$ such that $\R f_* C = 0$.
\end{itemize}

The rest of this section is devoted to the proof of the following statement.
\begin{thmm}\label{hearts}
	Assume to be working in Situation \ref{situ1} and assume in addition $f$ to have relative dimension at most one.
	Then the equivalence $\Phi$ between $D(Y)$ and $D(\orb{X})$ restricts to an equivalence of abelian categories between $\perv{0}{\Per}(Y/X)$ and $\Coh (\orb{X})$, while the equivalence $\otwit{\Phi}$ restricts to an equivalence between $\perv{-1}{\Per(Y/X)}$ and $\Coh(\orb{X})$.
\end{thmm}
\begin{rmk}
	Notice that the condition $\dim Y \times_X Y \leq n+1$ follows automatically from the condition on the fibres of $f$.
\end{rmk}
In particular $\perv{0}{\Per(Y/X)}$ is equivalent to $\perv{-1}{\Per(Y/X)}$.
We also point out that the composition $\otwit{\Phi}\twit{\Phi}^{-1}$ gives a non-trivial autoequivalence of $D(\orb{X})$, which seems related to the \emph{window shifts} of Donovan-Segal \cite{window}.
It might be worthwhile to compute this equivalence in explicit examples.

Let us now begin the proof of the theorem, which will be divided into small steps.
We start by considering $\twit{\Phi}$.
\smallskip
\stepz Given an object of the derived category, membership of either of the categories in question can be checked \'etale-locally on $X$ \cite[Proposition 3.1.6]{vdb}.
Thus, by base-changing over \'etale patches of $X$, we can reduce to the case where $X$ is affine and furthermore $\orb{X} = [V/G]$ with $V$ a smooth affine scheme and $G$ finite.
Moreover, the functors $\twit{\Phi}$ and $\twit{\Psi}$ (being Fourier-Mukai) commute with this base-change \cite[Proposition 6.1]{nahm}.

\smallskip
\stepz It suffices to prove $\twit{\Psi}(\Coh (\orb{X})) \subset \qerv{0}{\Per}(Y/X)$ because of the following well-known trick.
\begin{lem}
	Let $\curvy{A}$ and $\curvy{B}$ be two hearts relative to two bounded t-structures in a triangulated category.
	Then $\curvy{A} \subset \curvy{B}$ if and only if $\curvy{B} \subset \curvy{A}.$
\end{lem}
\begin{prf}
	Given an object $E$ let us denote by $H_\curvy{A}^i(E)$ (respectively $H_\curvy{B}^i(E)$) the $i$-th cohomology object relative to $\curvy{A}$ (resp.~$\curvy{B}$).
	Assume $\curvy{A} \subset \curvy{B}.$
	Let $E \in \curvy{B}$.
	As $E$ already lies in $\curvy{B}$ we have $E \simeq H^0_\curvy{B}(E)$ and $H^i_\curvy{B}(E) = 0$ for $i \neq 0.$
	Consider now the cohomology filtration of $E$ relative to $\curvy{A}$.
	As objects of $\curvy{A}$ are also in $\curvy{B}$, this filtration is also a filtration relative to $\curvy{B}$.
	By uniqueness of the cohomology objects we have $H^i_\curvy{A}(E) = H^i_\curvy{B}(E) = 0$ for $i\neq 0.$
	Thus, $E \in \curvy{A}$.
\end{prf}
\smallskip
\stepz To prove the mentioned inclusion we will exhibit two systems of generators (see definition below), one for $\qerv{0}{\Per}(Y/X)$ and one for $\Coh (\orb{X})$, and show that elements of the first system are sent to the second.
\begin{defn}
	Let $\cat{D}$ be a triangulated category and let $\curvy{A}$ be the heart of a bounded t-structure.
	A collection $\curly{P}$ of objects of $\curvy{A}$ is a \emph{system of projective generators} if, for all $A \in \curvy{A} \setminus \{0\}$ and all $P \in \curly{P}$, $\Ext_\cat{D}^\bullet(P,A)$ is concentrated in degree zero and for all $A \in \curvy{A}$ there exists $P_A \in \curly{P}$ such that $\Hom_\cat{D}(P_A,A) \neq 0$.
\end{defn}
By \cite[Lemma 3.2.4]{vdb}, when $X$ is affine, we have a system of generators $\curly{P}$ for $\qerv{0}{\Per}(Y/X)$ consisting of vector bundles $P$ such that
\begin{itemize}
	\item $\R^1f_*P = 0$,
	\item $P^\vee$ is generated by global sections.
\end{itemize}
For $\Coh (\orb{X})$ we also have a nice system of generators.
\begin{lem}
	The collection $\curly{Q}$ of vector bundles on $\orb{X}$ is a system of generators for $\Coh (\orb{X})$.
\end{lem}
\begin{prf}
	As we are working in the case $\orb{X} = [V/G]$, it is easy to reduce the problem to bundles on $V$.
	In fact, coherent sheaves on $\orb{X}$ are $G$-equivariant coherent sheaves on $V$.
	Given an equivariant vector bundle $P$ and an equivariant sheaf $E$ on $V$ we have that $G\text{-}\Ext^i_V(P,E) = \Ext_V^i(P,E)^G$, where the latter is the $G$-invariant part \cite[Section 4.1]{bkr}.
	As $V$ is affine, these groups vanish for $i>0$.
	
	Fix now an equivariant sheaf $E$, we want to find an equivariant vector bundle $P$ such that $\Hom_V(P,E)^G \neq 0$.
	By \cite[Lemma 4.1]{bkr} $\Hom_V(P,E)$ splits as a direct sum of $\Hom_V(P\otimes \rho, E)^G \otimes \rho$, where $\rho$ ranges among the irreducible representations of $G$.
	The claim thus follows as $P \otimes \rho$ is a vector bundle.
\end{prf}

\stepz We now conclude the proof by showing that elements of $\curly{P}$ are sent to elements of $\curly{Q}$.
First we remark that we can check whether a complex on $\orb{X} = [V/G]$ is a vector bundle by pulling back to the \'etale atlas $V \to [V/G]$.
Thus, if $P \in \curly{P}$, we are interested in the pullback of $\twit{\Phi}(P)$ to $V$.
This allows us to reduce to the setup of \cite{bkr}, where one has the following diagram.
\begin{center}
	\begin{tikzpicture}
		\matrix (m) [matrix of math nodes, row sep=3em, column sep=3em, text height=1.5ex, text depth=0.25ex]
		{&Z&\\
		&Y\times V&\\
		Y&&V\\
		&X&\\
		};
		\path[->,font=\scriptsize]
		(m-1-2) edge[bend right=20] node[auto,swap]{$p$}(m-3-1)
				edge[bend left=20] node[auto]{$q$}(m-3-3)
		(m-1-2) edge node[auto]{$i$} (m-2-2)
		(m-2-2) edge node[auto,swap]{$\pi_Y$} (m-3-1)
				edge node[auto]{$\pi_V$} (m-3-3)
		(m-3-1) edge node[auto,swap]{$f$}(m-4-2)
		(m-3-3) edge node[auto]{$g$}(m-4-2)
		;
	\end{tikzpicture}
\end{center}
Here $Z$ is the universal $G$-cluster for the action of $G$ on $V$, $q$ and $f$ are proper and birational, $p$ and $g$ are finite and $p$ is also flat.
Moreover, the quotient $\O_{Y \times V} \onto \O_Z$ is precisely the pullback, under the morphism $Y \times V \to Y \times [V/G] = Y \times \orb{X}$, of the universal quotient $\O_{Y \times \orb{X}} \onto \O_\orb{Z}$, which we used to define $\twit{\Phi}$.
It follows that applying $\twit{\Phi}$ followed by pulling back to $V$ is the same as applying $\R q_* p^*$.

We are thus reduced to checking that, given an element $P \in \curly{P}$, the complex $\R q_* p^* P$ is actually a vector bundle.
\begin{lem}
	Let $P \in \Coh(Y)$ satisfy $\R^1 f_* P = 0$.
	Then $\R q_* p^* P \in \Coh(\orb{X})$.
\end{lem}
\begin{prf}
	Notice that $\R q_* p^* P = \R \pi_{V,*} i_* p^* P = \R \pi_{V,*} (\pi_Y^* P \otimes \O_Z)$, where we made the standard identification $\O_Z = i_* \O_Z$.
	We point out that, as a consequence of our assumption on $f$, $\pi_{V,*}$ is of homological dimension at most one (we remind the reader that we work under the reduction done in the Step 1, in particular $X$ is affine).
	
	By tensoring the quotient $\O_{Y \times V} \onto \O_Z$ with $\pi_Y^*P$ we produce a surjection $\pi_Y^*P \onto \pi_Y^* P \otimes \O_Z$.
	Applying $\pi_{V,*}$ yields a surjection $\R^1\pi_{V,*} \pi_Y^* P \onto \R^1 \pi_{V,*} (\pi_Y^* P \otimes \O_Z)$.
	But $\R^1 \pi_{V,*} \pi_Y^* P = H^1(Y,P) \otimes_\C \O_V$ and $H^1(Y,P) = 0$ as $\R^1f_* P = 0$, hence the claim.
\end{prf}
\begin{lem}
	Let $P \in \curly{P}$, then $\R q_* p^* P$ is a vector bundle on $V$.
\end{lem}
\begin{prf}
	We know that the dual of $P$ is generated by global sections, hence there exists a short exact sequence
	\begin{align*}
		K \into \O^{\oplus m}_Y \onto P^\vee.
	\end{align*}
	From the fact that $P$ and $\O_Y$ are vector bundles it follows that $K$ is also a vector bundle.
	We therefore have a dual sequence
	\begin{align*}
		P \into \O^{\oplus m}_Y \onto K^\vee.
	\end{align*}
	It follows from the previous lemma, plus the fact that $q_* \O_Z = \O_V$, that applying $\R q_* p^*$ yields an exact sequence
	\begin{align*}
		q_* p^* P \into \O^{\oplus m}_V \onto q_* p^* K^\vee.
	\end{align*}
	To prove our claim it suffices to check that $\Ext^1_V(q_*p^*P,M) = 0$ for all modules $M$ on $V$.
	By the above short exact sequence this is the same as showing that $\Ext^2_V(q_*p^* K^\vee,M) = 0$ for all modules $M$.
	Using Grothendieck duality 
	for $q$ we have
	\begin{align*}
		\Ext_V^2(q_*p^* K^\vee,M)
		= \Ext_Z^2(p^*K^\vee,q^!M) 
		= H^2(Z,p^*K \otimes q^! M).
	\end{align*}
	The scheme $Z$ admits a finite and flat map to a smooth variety ($f: Z \to Y$) thus it is Cohen-Macaulay.
	Moreover, as $\dim Z - \dim V = 0$ and $q$ is of finite tor-dimension, the complex $q^!M$ is concentrated in non-positive degrees.
	As our assumption on $f$ implies that $H^i(Z,E) = 0$ for all $i > 1$ and all sheaves $E$, the hypercohomology spectral sequence tells us that $H^2(Z,p^*K \otimes q^!M) = 0$.
\end{prf}
The previous lemma concludes the first half of the proof.
As is often the case, the second half is much shorter than the first.
In fact, to prove the statement for $\perv{-1}{\Per(Y/X)}$ and $\Phi$, one need only notice the following:
\begin{itemize}
	\item $\otwit{\Phi} = \D \twit{\Phi} \D$, 
	\item the dual system $\curly{P}^\vee = \D \curly{P}$ is a system of generators for $\perv{-1}{\Per(Y/X)}$ \cite[3.2.3]{vdb},
	\item the system $\curly{Q}$ is self-dual $\D \curly{Q} = \curly{Q}$.
\end{itemize}
This concludes the proof and we can now move on to comparing the DT invariants of $\orb{X}$ and $Y$.
\begin{rmk}\label{o}
	For the next section, it will be important to know that $\Phi(\O_Y) = \O_\orb{X}$.
	We already know that $\Phi(\O_Y)$ is a vector bundle given by $\R q_*\O_\orb{Z}.$
	By restricting to the smooth locus of $X$ (viz.~to an open where $\Phi$ is the identity) we see that $\R q_* \O_\orb{Z}$ is in fact a line bundle.
	In turn this implies that $\Phi(\O_Y) = \O_\orb{X}$ as the unit $\O_\orb{X} \to \R q_* q^* \O_\orb{X}$ is an isomorphism.
	The same statement obviously holds for $\hat{\Phi} = \D \Phi \D$ as well.
\end{rmk}
%
%
%
\begin{rmk}\label{usefulcommute}
	It can be useful to know that when $Y$ and $X$ are projective the equivalences described above commute with pushing down to $X$.
	For example, let us check that $g_* \Phi = \R f_*$.
	We have $g_* \Phi = \R f_* \R p_* p^*$.
	If we proved that $\R p_* \O_\orb{Z} = \O_Y$, then by the projection formula we would be done.
	Thankfully, the previous remark together with Remark \ref{otherperversity} already tell us that $\R p_* \O_\orb{Z} = \Psi(\O_\orb{X}) = \O_Y$.
	%
	%
\end{rmk}


%% file: two.tex
\section{The Formula for DT Invariants} 
\label{sec:the_formula_for_dt_invariants}
We now impose further restrictions on our spaces.
\begin{situ}\label{situ2}
	Recall Situation \ref{situ1} and assume in addition $\orb{X}$ to be projective and of dimension three.
	Assume moreover $\orb{X}$ to be Calabi-Yau, i.e.~$\omega_\orb{X} \cong \O_\orb{X}$ and $H^1(\orb{X},\O_\orb{X})=0$.
	Finally, assume the crepant resolution $f\colon Y \to X$ of the previous section to have relative dimension at most one.
\end{situ}
\begin{rmk}
	We follow the convention where a Deligne-Mumford stack is projective if its coarse moduli space is.
	From the assumptions above it follows that $X$ is of dimension three, projective, Gorenstein with quotient singularities and with trivial canonical bundle.
	In turn it follows that $Y$ is Calabi-Yau of dimension three and that $X$ has rational singularities, and so $\R f_*\O_Y=\O_X$ \cite{kovacs}.
\end{rmk}
As the functor $\twit{\Phi}$ is more natural from the perspective of the McKay correspondence we shall focus on the zero perversity.
\begin{notation}
	We will drop the superscript $^{0}$ from $\perv{0}{\Per(Y/X)} =: \Per(Y/X)$.
\end{notation}

\addtocontents{toc}{\protect\setcounter{tocdepth}{0}}
\subsection{Reminder}\label{reminder}
Let us recall some definitions from the introduction.
We denote by $N(Y)$ the \emph{numerical K-group} of coherent sheaves of $Y$.
We remind ourselves that we can define a bilinear form on $K_0(\Coh (Y))$
\begin{align*}
	\chi(E,F) := \sum_k (-1)^k \dim_\C \Ext^k_Y(E,F)
\end{align*}
and that $N(Y)$ is obtained by quotienting out its radical.
Inside $N(Y)$ we can single out $F_1N(Y)$, which is the subgroup generated by sheaves supported in dimensions at most one.
We also define $F_{\exc}N(Y)$ to be the subgroup of $F_1N(Y)$ spanned by sheaves supported in dimension at most one and with derived pushforward to $X$ supported in dimension zero (see \cite[Section 4]{cala} for how this notion behaves well for perverse coherent sheaves).

To $Y$ one can also attach the numerical Chow groups $N_*(Y)$, which are the groups of cycles modulo numerical equivalence.
We write $N_{\leq 1}(Y) := N_1(Y) \oplus N_0(Y)$ and recall that $N_0(Y) \cong \Z$.
In \cite[Lemma 2.2]{tomcc} it is shown that the Chern character induces an isomorphism $F_1N(Y) \cong N_{\leq 1}(Y) \cong N_1(Y) \oplus \Z$, which allows us to pass from one group to the other.
Using this identification, $F_{\exc}N(Y)$ can be rewritten as
\begin{align*}
	F_{\exc} N(Y) = \left\{ (\beta,n)\in N_1(Y) \oplus \Z \st f_*\beta=0 \right\}
\end{align*}
where $f_*$ here stands for the proper pushforward on cycles (the subscript $\exc$ is short for \emph{exceptional}).

For the orbifold $\orb{X}$ we can also define a numerical K-group $N(\orb{X})$.
The functor $\Phi\colon D(Y) \to D(\X)$, with inverse $\Psi$, induces an isomorphism on the level of numerical K-groups.
\begin{align*}
	\phi: N(Y) \rightleftarrows N(\orb{X}): \psi
\end{align*}
The group $F_{\mr}N(\X)$ is defined to be $\phi(F_1N(Y))$ and we recall the diagram \eqref{delta} from the introduction, which expressed how all the numerical classes on $Y$ and $\X$ match up.
Before we may proceed, a technical remark is in order.
\begin{rmk}\label{sign issue crc}
	As mentioned in the introduction, we think of DT invariants as weighted Euler characteristics of the Hilbert scheme of a given Calabi-Yau threefold $M$, where the weight is given by Behrend's constructible function.
	The proof of the flop formula in \cite{cala} relies on the technology of motivic Hall algebras of Joyce.
	One of the technical points of this approach is that the Hilbert scheme $\Hilb(M)$ has a forgetful morphism $\sigma$ to the stack of coherent sheaves $\orb{M}$.
	Thus, on $\Hilb(M)$ there are two candidate constructible functions: the Behrend function $\nu_{\Hilb(M)}$ and the pullback $\mu = \sigma^* \nu_{\orb{M}}$ of the Behrend function on $\orb{M}$.
	Fortunately there is a simple relationship between the two, which unfortunately introduces some signs: namely if $\O_M \onto E$ is a quotient, with $E$ supported in dimension at most one then \cite[Theorem 3.1]{tomcc}
	\begin{align*}
		\nu_{\Hilb(M)}(\O_M \onto E) = (-1)^{\chi(E)} \mu(\O_M \onto E) = (-1)^{\chi(E)} \nu_{\orb{M}} (E)
	\end{align*}
	where $\chi(E) = \chi(\O_M,E)$ is the Euler characteristic of the sheaf $E$.
	Moreover if $f\colon U \to Z$ is an open immersion then $\nu_U = f^*\nu_Z$.
	This is relevant for the ``partial'' invariants.
\end{rmk}
The \emph{DT number} of $Y$ of class $\alpha \in F_1N(Y)$ is given by 
\begin{align*}
	\DT_Y(\alpha) := \chi_{\topp}\left( \Hilb_Y(\alpha), \nu \right)
\end{align*}
but for convenience, we give a name to the numbers obtained by weighing with $\mu$ as well, namely
\begin{align*}
	\correct{\DT}_Y(\alpha) := \chi_{\topp}\left( \Hilb_Y(\alpha), \mu \right) = (-1)^{\chi(\alpha)} \DT_Y(\alpha).
\end{align*}
In the introduction we also mentioned that we package all these numbers into a generating series
\begin{align*}
	\DT(Y) := \sum_{\alpha \in F_1N(Y)} \DT_Y(\alpha)q^\alpha
\end{align*}
and similarly for $DT_0(Y)$, $DT_{\exc}(Y)$ and the corresponding underlined versions.
Finally, we also define
\begin{align*}
	DT_{\exc}^\vee(Y) := \sum_{\substack{(\beta,n) \in N_1(Y)\oplus \Z \\ f_*\beta = 0}} DT_Y(-\beta,n) q^{(\beta,n)}
\end{align*}
and the corresponding $\correct{DT}^\vee_{\exc}(Y)$.

Recall now the category $\Per(Y/X)$ of perverse coherent sheaves from the previous section.
The structure sheaf $\O_Y$ belongs to $\Per(Y/X)$ and there is a moduli space $\PHilb(Y/X)$ parameterising quotients of $\O_Y$ in $\Per(Y/X)$ \cite[Section 6]{tomflops}.
This space splits into open and closed components $\PHilb_{Y/X}(\alpha)$, for each numerical class $\alpha$, parameterising quotients $\O_Y \onto P$, with $[P] = \alpha$.

We can define a \emph{perverse} DT number of $Y$ over $X$ of class $\alpha$ as the weighted Euler characteristic
\begin{align*}
	\correct{\DT}_{Y/X}(\alpha) := \chi_{\topp}\left( \PHilb_{Y/X}(\alpha), \mu \right)
\end{align*}
where $\mu$ is the pullback of the Behrend function of the stack of perverse coherent sheaves on $Y$.
We also collect these numbers into a generating series
\begin{align*}
	\correct{\DT}(Y/X) := \sum_{\alpha \in F_1N(Y)} \correct{\DT}_{Y/X}(\alpha) q^\alpha, \quad \quad
	\correct{\DT}_{\exc}(Y/X) := \sum_{\alpha \in F_{\exc}N(Y)} \correct{\DT}_{Y/X}(\alpha) q^\alpha.
\end{align*}

The orbifold side brings no surprises.
We once again define DT numbers by taking weighted Euler characteristics and gather them in a generating series
\begin{align*}
	\DT_{\orb{X}}(\alpha) := \chi_{\topp}\left( \Hilb_\orb{X}(\alpha), \nu \right),
	\quad \quad
	\DT_{\mr}(\orb{X}) := \sum_{\alpha \in F_{\mr}N(\orb{X})} \DT_{\orb{X}}(\alpha) q^\alpha,
	\quad \quad
	\DT_0(\X) := \sum_{\alpha \in F_0N(\X)} DT_\X(\alpha) q^\alpha.
\end{align*}
\begin{rmk}\label{sign issue again crc}
	The analogue of Remark \ref{sign issue crc} for $\orb{X}$ still holds.
	That is, the following identity holds
	\begin{align*}
		\chi_\text{top} \left( \Hilb_{\orb{X}}(\alpha), \mu \right) = (-1)^{\chi(\alpha)} \chi_\text{top}\left( \Hilb_{\orb{X}}(\alpha), \nu \right)
	\end{align*}
	where $\mu$ is the pullback of the Behrend function of the stack of coherent sheaves on $\orb{X}$.
	To prove this, one can choose an appropriate divisor $D$ on the coarse space $X$, and its pullback to $\orb{X}$ plays the role of $H$ in the proof of \cite[Theorem 3.1]{tomcc}.
	The affine $U$ can then be chosen to be an \'{e}tale open in $\orb{X}$, so that \cite[Lemma 3.2]{tomcc} can be applied.
\end{rmk}
Because of this remark, we can define the underlined version of $\DT_{\mr}(\orb{X})$ and the identity above translates to
\begin{align*}
	\correct{\DT}_{\orb{X}}(\alpha) = (-1)^{\chi(\alpha)}\DT_\orb{X}(\alpha).
\end{align*}

\begin{rmk}\label{quot}
	Form the previous section we know that the Fourier-Mukai equivalences $\Phi$ and $\Psi$ restrict to an equivalence of abelian categories between $\Per(Y/X)$ and $\Coh (\orb{X})$.
	Using Remark \ref{o}, which tells us that $\Phi(\O_Y) = \O_\orb{X}$, we have an induced isomorphism of Quot functors (or Hilbert schemes), hence
	\begin{align*}
		\Hilb_\orb{X}(\alpha) \simeq \PHilb_{Y/X}(\psi(\alpha)).
	\end{align*}
\end{rmk}
\begin{rmk}\label{technical}
	Before we state the theorem, we point out a technical detail.
	In \cite{tomcc} the generating series $\DT(Y)$ is interpreted as belonging to an algebra $\C[\Delta]_\Phi$ (where $\Delta \subset F_1N(Y)$ is the positive cone of classes $[E]$ with $E \in \Coh(Y)$),  whose elements consist of formal series
	\begin{align*}
		\sum_{(\beta,n) \in \Delta \subset N_1(Y) \oplus \Z} a_{(\beta,n)} q^{(\beta,n)}
	\end{align*}
	where the $a_{(\beta,n)}$ are complex coefficients such that, for a fixed $\beta$, $a_{(\beta,n)}=0$ for $n$ very negative.
	
	A similar interpretation is given in \cite{cala} for the generating series $\correct{\DT}(Y/X)$, which now belongs to an algebra\footnote{Here the subscripts $\Phi$ and $\Lambda$ are just notation and stand for entirely parallel constraints. Also, to be pedantic, in \cite{cala} $\mathbb{Q}$ was used in place of $\mathbb{C}$. However the latter is obtained by the former by tensoring with $\mathbb{C}$.} $\C[\pepe{\Delta}]_\Lambda$ (where $\pepe{\Delta} \subset F_1N(Y)$ is the positive cone of classes $[P]$ with $P \in \Per(Y/X)$), which is the analogous of $\C[\Delta]_\Phi$ for $\Per(Y/X)$.
	The generating series $\correct{\DT}(Y)$ can also be seen as an element of $\C[\pepe{\Delta}]_\Lambda$ and all the identities we write down below should be interpreted as taking place within this algebra.
\end{rmk}
In light of \ref{quot} our main theorem is now immediate.
\begin{thmm}\label{main}
	Assume to be working in Situation \ref{situ2}.
	For each $\alpha \in F_1N(Y)$ we have
	\begin{align*}
			\correct{\DT}_{Y/X}(\alpha) = \correct{\DT}_\orb{X}(\phi(\alpha)).
	\end{align*}
	In particular, the following formulae hold.
	\begin{align*}
		\correct{\DT}_{\mr}(\orb{X}) = \correct{\DT}(Y/X), \quad \quad
		\correct{DT}_0(\X) = \correct{DT}_{\exc}(Y/X)
	\end{align*}
	after an identification of variables via $\phi$.
\end{thmm}
As one can plainly see, up to this point we never had to restrict to ``partial'' DT invariants or to zero-dimensional sheaves on $\X$.
The crepant resolution conjecture truly boils down to the study of these exotic invariants $\correct{DT}(Y/X)$.
\begin{cor}\label{mr=exc}
	Assume to be working in Situation \ref{situ2} and recall the identification of variables from the previous theorem.
	The following formula is true.
	\begin{align}\label{onceagain}
		\DT_0(\orb{X}) = \frac{\DT^\vee_{\exc}(Y)\DT_{\exc}(Y)}{\DT_0(Y)}
	\end{align}
\end{cor}
\begin{prf}
	First we notice that we can get rid of the underlines thanks to Remarks \ref{sign issue crc} and \ref{sign issue again crc}.
	From the previous theorem the statement we wish to prove is equivalent to proving that $\correct{\DT}_{\exc}(Y/X)$ is equal to the right hand side of \eqref{onceagain} (modulo the underlines).
	To apply \eqref{flopfo} from the introduction, we need to check that the assumptions of \cite[Remark 1.9]{cala} are satisfied.
	
	This is a simple issue of perversities ($p=0$, the one we use, versus $p=-1$, the one dominating \cite{cala}).
	The only thing to prove is that the stack of perverse coherent sheaves is locally isomorphic to the stack of coherent sheaves (via an \emph{exact} functor).
	However, using the Fourier-Mukai equivalence $\otwit{\Phi}$, we have that the stack parameterising objects in $\perv{0}{\Per(Y/X)}$ is isomorphic to the stack parameterising objects in $\perv{-1}{\Per(Y/X)}$.
	As $\otwit{\Phi}$ is also an exact functor, all the constructions of \cite{cala} go through and \eqref{flopfo} does indeed hold.
\end{prf}

\medskip
If we knew \eqref{c8} to be true, then Theorem \ref{main} would immediately imply \eqref{c0}.
Combined with Corollary \ref{mr=exc} we would derive a proof of \eqref{c1} \cite[Conjecture 1]{crc}.

Nonetheless, an unconditional variant can be proved \eqref{c0del}.
Let $\PHilb^\partial \subset \PHilb$ denote the open subspace consisting of those epimorphisms $\O_Y \to E$, where $\dim \supp E \leq 1$.
We then have a corresponding DT series
\begin{align*}
	\correct{DT}^\partial(Y/X) = \sum_{\alpha \in F_1N(Y)} \correct{DT}^\partial_{Y/X}(\alpha) q^\alpha,
	\quad \quad \quad
	\correct{DT}^\partial_{Y/X}(\alpha) := \chi_{\top}\left( \PHilb_{Y/X}^\partial(\alpha) ,\mu\right).
\end{align*}
Over $\X$, we let $\Hilb^\partial_\X(\alpha) \subset \Hilb_\X(\alpha)$ the image of $\PHilb_{Y/X}^\partial(\alpha)$ under $\Phi$.
The associated DT series is of course
\begin{align*}
	\correct{DT}_{\mr}^\partial(\X) = \sum_{\alpha \in F_{\mr}N(\X)} \correct{DT}^\partial_{\X}(\alpha) q^\alpha,
	\quad \quad \quad
	\correct{DT}^\partial_{\X}(\alpha) := \chi_{\top}\left( \Hilb_{\X}^\partial(\alpha) ,\mu\right).
\end{align*}
Finally, let $\Hilb^\partial_{Y}(\alpha) \subset \Hilb_Y(\alpha)$ be the subspace of $\Hilb_Y(\alpha)$ parameterising those quotients $\O_Y \to Q$, having perverse cokernel in $\Coh_{\leq 1}[1]$ (see \cite[Remark 3.5]{cala}).
For the last time, we have the generating series that comes with this moduli space.
\begin{align*}
	\correct{DT}^\partial(Y) = \sum_{\alpha \in F_1N(Y)} \correct{DT}^\partial_{Y}(\alpha) q^\alpha,
	\quad \quad \quad
	\correct{DT}^\partial_{Y}(\alpha) := \chi_{\top}\left( \Hilb_{Y}^\partial(\alpha) ,\mu\right).
\end{align*}

\medskip
Theorem \ref{main} specialises to the following.
\begin{thmm}
	Assume to be working in Situation \ref{situ2}.
	For each $\alpha \in F_1N(Y)$ we have
	\begin{align*}
			\correct{\DT}^\partial_{Y/X}(\alpha) = \correct{\DT}^\partial_\orb{X}(\phi(\alpha)).
	\end{align*}
	In particular, the following formula holds (after an identification of variables via $\phi$).
	\begin{align*}
		\correct{\DT}^\partial_{\mr}(\orb{X}) = \correct{\DT}^\partial(Y/X)
	\end{align*}
\end{thmm}
Using \cite[Theorem 3.30]{cala}, we have our last formula.
\begin{cor}\label{corollario}
	Assume to be working in Situation \ref{situ2} and recall the identification of variables using $\phi$.
	Then
	\begin{align*}
		DT^\partial_{\mr}(\X) = \frac{ DT^\partial(Y) DT^\vee_{\exc}(Y) } {DT_0(Y)},
		\quad \quad \quad
		\frac { DT^\partial_{\mr}(\X) }{DT_0(\X)} = 
		\frac { DT^\partial(Y) }{ DT_{\exc}(Y) }.
	\end{align*}
\end{cor}